\documentclass[graybox]{svmult}

\usepackage{type1cm}        
\usepackage{makeidx}         
\usepackage{graphicx}        
\usepackage{multicol}        
\usepackage[bottom]{footmisc}

\usepackage{subcaption}

\usepackage{newtxtext}       %
\usepackage[varvw]{newtxmath}       

\usepackage{orcidlink}

\makeindex             

\begin{document}

\title*{Model selection focusing on longtime behavior of differential equations}
 \titlerunning{Model selection using longtime behavior}
\author{Cordula Reisch\orcidlink{0000-0003-1442-1474} and Hannah Burmester\orcidlink{0009-0003-9802-6902}}
\institute{Cordula Reisch, Hannah Burmester \at TU Braunschweig, Institute for Partial Differential Equations,\email{c.reisch@tu-braunschweig.de}}

\maketitle

\abstract*{
Modeling biological processes is a highly demanding task because not all processes are fully understood. 
Mathematical models allow us to test hypotheses about possible mechanisms of biological processes. The mathematical mechanisms oftentimes abstract from the biological micro-scale mechanisms. 
Experimental parameter calibration is extremely challenging as the connection between abstract and micro-scale mechanisms is unknown. Even if some microscopic parameters can be determined by isolated experiments, the connection to the abstract mathematical model is challenging.
We present ideas for overcoming these difficulties by using longtime characteristics of solutions for, first, finding abstract mechanisms covering large-scale observations and, second, determining parameter values for the abstract mechanisms. 
The parameter values are not directly connected to experimental data but serve as a link between known mechanisms and observations.
The framework combines machine learning techniques with the characteristic solution behavior of differential equations. 
This setting gives insight into challenges by using rare data only that can later be used for partial differential equations.
}

\abstract{Modeling biological processes is a highly demanding task because not all processes are fully understood. 
Mathematical models allow us to test hypotheses about possible mechanisms of biological processes. The mathematical mechanisms oftentimes abstract from the biological micro-scale mechanisms. 
Experimental parameter calibration is extremely challenging as the connection between abstract and micro-scale mechanisms is unknown. Even if some microscopic parameters can be determined by isolated experiments, the connection to the abstract mathematical model is challenging.
We present ideas for overcoming these difficulties by using longtime characteristics of solutions for, first, finding abstract mechanisms covering large-scale observations and, second, determining parameter values for the abstract mechanisms. 
The parameter values are not directly connected to experimental data but serve as a link between known mechanisms and observations.
The framework combines machine learning techniques with the characteristic solution behavior of differential equations. 
This setting gives insight into challenges by using rare data only that can later be used for partial differential equations.
}

\section{Motivation and problem formulation}

Natural processes in life sciences generally involve many cells, substances, and particles interacting in a non-linear complex behavior. Understanding the whole system in detail is a very challenging task that is often tackled by first studying the interplay of fewer selected agents and their interactions. Even this study of a subset of all interactions is difficult as the mechanisms in an isolated system may be different and experiments are expensive and hard to realise due to the small length scales of the agents. 
At the same time, some observations on the system behavior take place on very different length scales and time scales and may be observed by measurements of large-scale quantities. The connection between the small-scale mechanisms and the large-scale outcome is usually non-monotonous and non-linear. 
Observations on the large scale can rather be interpreted as qualitative data, providing information on the general behavior of the system, than as quantitative data that gives input-output information on the underlying small-scale system. 

An exemplary problem for the described situation is viral liver inflammation like hepatitis B. The reasons for persisting or chronic liver inflammation are still not fully understood, \cite{Thomas2016}. Meanwhile, the interactions of various immune cells, the liver tissue, and the virus are complex. Time- and space-resolved data is not available because recordings of living human tissue are not possible. Examples of large-scale data are inflammatory values of the blood or, in an even larger setting, the development of chronic infections lasting for months. 
Mathematical models for hepatitis B, e.g. \cite{reisch-chronification,reisch-chemotactic}, provide a possibility to test different hypotheses of involved mechanisms. By taking the large-scale observations as qualitative data, various models with different mechanisms can be compared in their feasibility to reproduce this data. 

Here, we follow a top-down approach for model selection, starting with the qualitative data and giving a selection of possible models with the same qualitative solution behavior as an output. The model selection is carried out with machine learning approaches, following the idea of physics-informed neural networks (PINNs). 
We propose an algorithm for selecting mechanisms out of a class of mechanisms for exemplary ordinary differential equation systems, extending the ideas in \cite{reisch-automative}.

\section{Idea of model selection}

There are several approaches for learning (parts of) differential equations, focusing on different aspects. 
While the most common approaches as described in \cite{BruntonKutz-book} or in \cite{Meidani_2021,Rackauckas} use quite a lot of data, other approaches like in \cite{Schmiester_qualitative} and
\cite{Schmiester_parametrization} start to develop methods for less or qualitative data. 
Nevertheless, \cite{Alber} highlights the difficulties for machine learning approaches dealing with systems from biology, in particular acting on various scales. 
We use the technical methods in \cite{SCIANN} for identifying parameters of the differential equation system
\begin{align}
\begin{aligned}\label{eq:competitionmodel}
    u' &= u (1- a_1 u - a_2 v), \quad \quad
     &v'= r v ( 1- b_1 u - b_2 v)
\end{aligned}
\end{align}
describing the competition of species $u$ and $v$ with shared resources. Initial values are $u(0)=u_0>0$ and $v(0)=v_0>0$. Depending on the parameters $a_1, a_2, b_1, b_2, r$, the system tends towards a co-existence steady state, where ${(u^\star, v^\star) \neq (0,0)}$, or a single-species steady state with either $u^\star =0$ or $v^\star =0$.
The occurrence of one of the states, coexistence or distinction of one species, is qualitative data for this life science system. In contrast, quantitative data is the time-dependent number of species for many time steps over the whole dynamical process. Here, we interpret time-dependent data for single time steps during the quasi-stationary phase of the dynamics as well as quantitative data.

 \begin{figure}
 \begin{subfigure}[t]{.45\linewidth}
    \centering\includegraphics[width=\linewidth]{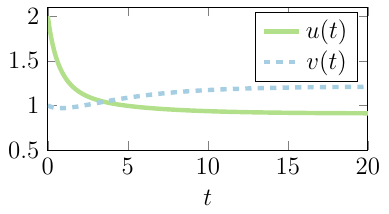}
    \caption{
    Coexistence:
    ${a_1=0.7}$, ${a_2=0.3}$, ${b_1=0.3}$, ${b_2=0.6}$.}
    \label{fig:coexistence}
  \end{subfigure}
  \hfill
  \begin{subfigure}[t]{.45\linewidth}
    \centering\includegraphics[width=\linewidth]{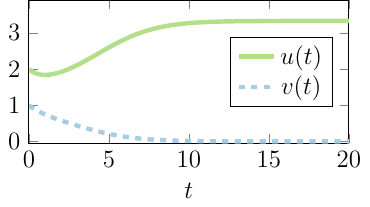}
    \caption{
    Survival of one species: 
    ${a_1=0.3}$, ${a_2=0.6}$, ${b_1=0.7}$, ${b_2=0.3}$.}
    \label{fig:distinction}
  \end{subfigure}
  \caption{Solutions of the competition model for different parameter settings. Initial values are $(u_0,v_0)=(2,1)$ and $r=0.5$ is fixed.}
     \label{fig:examples_competition}
 \end{figure}

Learning equations by machine learning techniques can be handled in multiple ways. Here, we focus on learning equations, not solutions, which also can be done with neural networks. 
Packages like SciANN, \cite{SCIANN}, allow to learn parameters in an equation like \eqref{eq:competitionmodel}. On top, packages like PySINDy, \cite{PySindy}, have already been implemented to favor sparse models by including the number of included mechanisms in the loss function.  

We combine the two approaches by using the simple architecture of SciANN and reduce the number of included mechanisms step-by-step. 
Therefore, we start with a model that includes more mechanisms, namely
\begin{align} \label{eq:competition_plus}
\begin{aligned}
    u' &= u (1- a_1 u - a_2 v) + \alpha_0 +\alpha_1 u +\alpha_2 v + \alpha_3 u^2v + \alpha_4 uv^2,\\
    v'&= r v ( 1- b_1 u - b_2 v) +\beta_0 +\beta_1 u +\beta_2 v + \beta_3 u^2 v+ \beta_4 uv^2.
\end{aligned}
\end{align}
Numerical simulations of system \eqref{eq:competitionmodel} provide time-resolved data which is used for training, testing, and evaluation of the learned model based in \eqref{eq:competition_plus}. 

The technical setting of SciANN is listed in Table~\ref{tab:technical}.

\begin{table}
    \centering
        \caption{Parameter setting for the neural networks in SciANN. All layers had the same number of neurons. }
    \label{tab:technical}
    \begin{tabular}{p{2cm}p{1cm}|p{0.1cm}p{2.5cm}p{1cm}|p{0.1cm}p{3cm}p{1.5cm}}
    \hline\noalign{\smallskip}
                layers & 5 & &batch size & 100 && activation function & tanh  \\
       neurons & 10 & &epochs & 5000 & &loss function &  MSE \\
       & & & number data points & 100& &optimization algorithm & Adams\\
    \end{tabular}
\end{table}

The model selection follows the procedure
\begin{enumerate}
    \item Fit the (positive) parameter in the model to the given data. 
    \item Delete the mechanism with the smallest parameter.
     \item Define the new model without the deleted mechanisms
    \item[$\circ$] Repeat 1.-3. as long as the number of selected mechanisms is above a threshold $\varepsilon$ of included mechanisms.  
\end{enumerate}

Following this algorithm, a model is selected step-by-step from the family of mechanisms given in system \eqref{eq:competition_plus}.
The used threshold is $\varepsilon=6$.

\section{Learning ordinary differential equations with a focus on longtime behavior}

The algorithm was tested for two scenarios where the time interval for the data points was restricted differently. 

Before we apply the selection algorithm, we start with learning the parameters given in Eq.~\eqref{eq:competitionmodel} directly. Afterwards, we regard the same data scenarios for the step-by-step learning following the proposed algorithm. 

\subsection{Learning the parameters of a given model}\label{sec:31}

\begin{figure}
 \begin{subfigure}[t]{.45\linewidth}
    \centering\includegraphics[width=\linewidth]{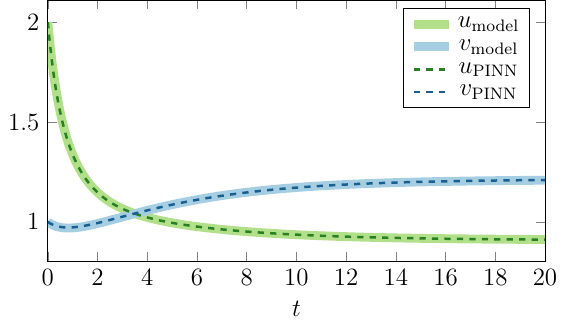}
    \caption{Data points $u_\mathrm{data}, v_\mathrm{data}$ from the whole time interval $[0,20]$.}
     \label{fig:PINN_full}
  \end{subfigure}
  \hfill
  \begin{subfigure}[t]{.45\linewidth}
    \centering\includegraphics[width=\linewidth]{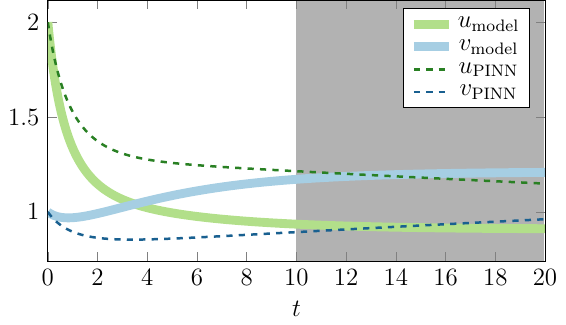}
    \caption{Data points $u_\mathrm{data}, v_\mathrm{data}$ from the time interval $[10,20]$, highlighted in grey.}
       \label{fig:PINN_stat}
  \end{subfigure}
  \caption{Solution $u_\mathrm{model}, v_\mathrm{model}$ of model~\eqref{eq:competitionmodel} and solution $u_\mathrm{PINN}, v_\mathrm{PINN}$ of the model with the learned parameters in  Table~\ref{tab:selection_fulldomain} and Table~\ref{tab:selection_latedomain}, respectively.}
     \label{fig:comparison_PINN}
 \end{figure}

The first scenario is the availability of data in the whole time interval. We compare the results by use of the mean squared error of the solutions. For visualization, we plot in Fig.~\ref{fig:PINN_full} the solution $u_\mathrm{PINN}$ of the ordinary differential equation with the learned parameters and the solution $u_\mathrm{model}$ from which we extracted the data 
\begin{equation}
    u_\mathrm{data}= \{ (u_0, v_0)^\mathrm{T} \} \cup \{ (u_\mathrm{model}(t_i), v_\mathrm{model}(t_i))^\mathrm{T} : i=1,\dots, N \}.
\end{equation}
The difference between the two scenarios is the  interval for $t_i$, which is $[0,20]$ for the whole time interval and $[10,20]$ for the time interval with a quasi-stationary solution. 
In both scenarios, the initial data is given as an extra information.

The availability of qualitative data is interpreted as having data only from the time interval where the solution is already quasi-stationary, see Fig.~\ref{fig:PINN_stat}. 

The solutions of the PINN approach in Fig.~\ref{fig:comparison_PINN} show a good match between the solution providing the data and the solution with the learned parameters. 
The data points in Fig.~\ref{fig:PINN_stat} are all from the quasi-stationary time interval. The difference between the model solution and the learned solution with the PINN approach is larger for small time $t$. The quality of the solution with data from the whole time interval is higher in total. 
The error for the learned solution is given in Fig.~\ref{fig:MSE_selection}.

\subsection{Applying the model selection algorithm}

We start again with the first scenario with data from the whole time interval, and reduce the number of mechanisms in every step following the selection algorithm. Tab.~\ref{tab:selection_fulldomain} shows the parameters in every selection step. 
We stop the selection algorithm when we are having $\varepsilon=6$ parameters left, which is one additional parameter compared to the model in Eq.~\eqref{eq:competitionmodel}.

\begin{table}
    \centering
        \caption{Learned parameters of Eq.~\eqref{eq:competition_plus} with data from the whole time interval.  }
        \label{tab:selection_fulldomain}
\begin{tabular}{c | lllllllll | l | r}
& &&&&step &&&&\\
    parameter $\quad$   & 1     & 2     & 3     & 4     & 5     & 6     & 7     & 8     & 9  & Sec. \ref{sec:31} &$\ $ data\\ \hline
 $r$   & 0.491 $\ $ & 0.443$\ $ & 0.471 $\ $& 0.489 $\ $& 0.668$\ $ & 0.683 $\ $& 0.516 $\ $& 0.565$\ $ & {\bf 0.500} $\ $ & 0.478 & 0.5\\
 $a_1$ & 0.731 & 0.764 & 0.770 & 0.771 & 0.842 & 0.732 & 0.718 & 0.625 &  {\bf 0.623} & 0.699 & 0.7 \\
 $a_2$ & 0.789 & 0.912 & 0.902 & 0.884 & 0.788 & 0.772 & 0.773 & 0.687 &  {\bf 0.681 }& 0.301 & 0.3  \\
 $b_1$ & 0.914 & 0.898 & 0.866 & 0.856 & 0.788 & 0.746 & 0.725 & 0.679 &  {\bf 0.289} & 0.295 & 0.3\\
 $b_2$ & 1.150 & 0.935 & 0.927 & 0.901 & 0.833 & 0.829 & 0.748 & 0.709 &  {\bf 0.609} & 0.603 & 0.6 \\
 $\alpha_0$ & 0.001 & 0.245 & 0.232 & 0.242 & 0.284 & 0.310 & 0.299 & 0.357 &  {\bf 0.349 }& & 0\\
 $\beta_0$ & 0.002 & 0.248 & 0.252 & 0.277 & 0.291 & 0.282 & 0.343 & 0.317 & / & & 0   \\
 $\alpha_1$ & 0.228 & 0.187 & 0.166 & 0.168 & 0.211 & 0.251 & 0.250 &  /    & /   & & 0  \\
 $\beta_1$ & -0.   & /     & /     & /     &  /    & /     & /     & /     & /    & & 0 \\
 $\alpha_2$ & 5.7e-5& 0.120 & 0.127 & 0.115 & /     & /     & /     & /     & /    & & 0 \\
 $\beta_2$ & 0.063 & 0.150 & 0.172 & 0.192 & 0.238 & 0.223 & /     & /     & /    & & 0  \\
 $\alpha_3$ & 0.062 & 0.162 & 0.180 & 0.165 & 0.174 & /     & /     & /     & /    & & 0 \\
 $\beta_3$ & -0.   & /     & /     & /     & /     & /     & /     & /     & /    & & 0 \\
 $\alpha_4$ & 0.222 & 2.9e-4& /     & /     & /     & /     & /     & /     & /   & & 0  \\
 $\beta_4$ & 0.486 & 0.056 & 0.047 & /     & /     & /     & /     & /     & /    & & 0 
 \end{tabular}
 \end{table}

\begin{figure}
 \begin{subfigure}[t]{.45\linewidth}
    \centering\includegraphics[width=\linewidth]{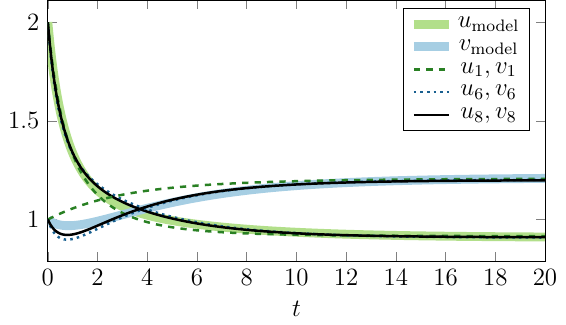}
    \caption{Data points $u_\mathrm{data}, v_\mathrm{data}$ from the whole time interval $[0,20]$.}
     \label{fig:selection_fulldomain}
  \end{subfigure}
  \hfill
  \begin{subfigure}[t]{.45\linewidth}
    \centering\includegraphics[width=\linewidth]{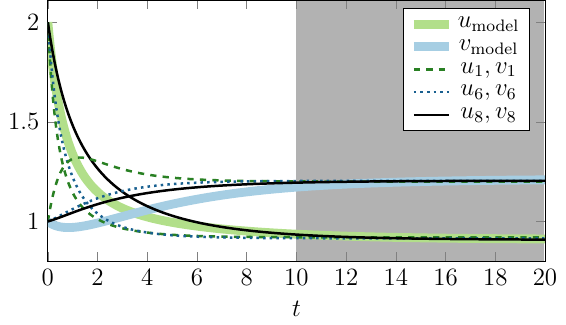}
    \caption{Data points $u_\mathrm{data}, v_\mathrm{data}$ from the time interval $[10,20]$, highlighted in grey.}
       \label{fig:selection_latedomain}
  \end{subfigure}
  \caption{Solution $u_\mathrm{model}, v_\mathrm{model}$ of model~\eqref{eq:competitionmodel} and solutions $u_\mathrm{PINN}, v_\mathrm{PINN}$ of the step-by-step selected models with parameters in Table~\ref{tab:selection_fulldomain} and Table~\ref{tab:selection_latedomain}, respectively. }
     \label{fig:comparison_selection}
 \end{figure}

The selection process in Tab.~\ref{tab:selection_fulldomain} shows that all included mechanisms are learned and additionally the parameter $\alpha_0$. 
The learned model reads
\begin{align} \label{eq:learned_full}
\begin{aligned}
    u' &= u (1- a_1 u - a_2 v) + \alpha_0 , \quad \quad 
    & v'= r v ( 1- b_1 u - b_2 v) ,
\end{aligned}
\end{align}
and the largest difference between the data-providing parameter and the learned parameter is for $a_2$. 
Fig.~\ref{fig:selection_fulldomain} shows the solutions of the models in different selection steps. 
The model solutions $v_1$ and $v_6$ in the first and sixth selection steps are worse than the solution $v_9$ of the reduced model in the final selection step. 
This result is surprising from an approximation point of view. In the first step of the selection process, the number of free parameters is the largest and therefore the fitting should be best. As the setting of the learning process is the same for all selection steps, a reduced number of parameters can be an advantage for fitting those parameters best. 

In comparison to the approach with fully available data, we take again the data points only from the time interval where the solution in Fig.~\ref{fig:coexistence} is quasi-stationary. Tab.~\ref{tab:selection_latedomain} shows the parameters in every selection step.

\begin{table}
    \centering
        \caption{Learned parameters of Eq.~\eqref{eq:competition_plus} with data from time interval $[10,20]$.  }
        \label{tab:selection_latedomain}
\begin{tabular}{c | llllllll | l  | r}
para-& &&&step &&&&& &\\
    meter& 1   & 2     & 3     & 4     & 5     & 6     & 7     & 8     & Sec. \ref{sec:31}  &  data \\ \hline
 $r$   & 0.788$\ $ & 0.630 $\ $ & 0.406 $\ $& 0.512$\ $ & 0.573$\ $ & 0.598 $\ $& 0.634$\ $ & {\bf 0.677} $\ $& 0.475 $\ $ & 0.5 \\
 $a_1$ & 1.207 & 0.625 & 0.623 & 0.851 & 0.790 & 0.771 & 0.653 & {\bf 0.488} & 0.487  & 0.7 \\
 $a_2$ & 0.643 & 0.872 & 0.888 & 0.797 & 0.770 & 0.757 & 0.680 & {\bf 0.463} & 0.463  & 0.3  \\
 $b_1$ & 0.585 & 0.610 & 0.739 & 0.819 & 0.805 & 0.683 & 0.672 & {\bf 0.664} & 0.456  & 0.3 \\
 $b_2$ & 1.399 & 1.195 & 0.916 & 0.753 & 0.764 & 0.659 & 0.678 & {\bf 0.652} & 0.480  & 0.6  \\
 $\alpha_0$ & 0.007& 0     & /     & /     & /     & /     & /     & /     & /     & 0\\
 $\beta_0$ & 0     & /     & /     & /     & /     & /     & /     & /     & /     & 0 \\
 $\alpha_1$ &0.258 & 0.251 & 0.291 & 0.218 & 0.250 & 0.245 & /     & /     & /     & 0    \\
 $\beta_1$ & 0.449 & 0.402 & 0.344 & 0.194 & 0.222 & 0.329 & 0.359 & {\bf 0.348} & /     & 0  \\
 $\alpha_2$ & 0.007& 0.0004& 0.045 & 0.280 & 0.305 & 0.283 & 0.316 & /     & /     & 0  \\
 $\beta_2$ & 0     & /     &   /   &   /   & /     & /     & /     & /     & /     & 0   \\
 $\alpha_3$ & 0.058& 0.329 & 0.253 & 0.137 & /     & /     & /     & /     & /     & 0 \\
 $\beta_3$ & 0     & /     & /     & /     & /     & /     & /     & /     & /     & 0  \\
 $\alpha_4$ & 0.376& 0.001 & 0.003 & /     & /     & /     & /     & /     & /     & 0\\
 $\beta_4$ & 0.557 & 0.289 & 0.049 & 0.172 & 0.189 & /     & /     & /     & /     & 0   
 \end{tabular}
 \end{table}

\begin{figure}
\sidecaption
    \includegraphics[width=0.6 \textwidth]{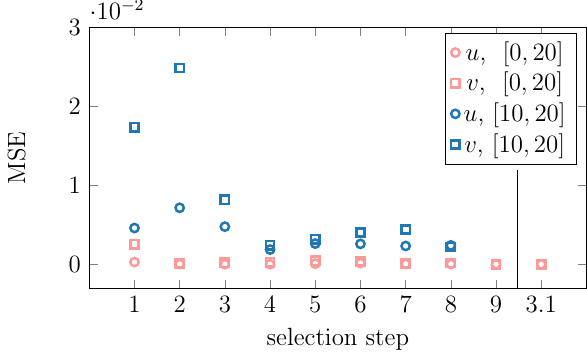}
    \caption{Mean squared error for models during the selection step. The number of selection steps differs due to the elimination of multiple zeros in one step. The last column gives the MSE for the parameter fitting without selection in Sec.~\ref{sec:31}. }
    \label{fig:MSE_selection}
\end{figure}

By using only data from the quasi-stationary phase of the dynamics, the selected model is given by 
\begin{align} \label{eq:selected_2}
\begin{aligned}
    u' &= u (1- a_1 u - a_2 v)  , \quad \quad 
    &v'= r v ( 1- b_1 u - b_2 v) + \beta_1 u.
\end{aligned}
\end{align}
The parameters differ without a clear tendency from the data parameters.  
Fig.~\ref{fig:selection_latedomain} shows the solutions of the models in different selection steps. 
The data points are all located in the grey region, and the approximation of the solutions in this domain is quite accurate. 
The solution in the last selection step approximates the dynamical behavior in the beginning better than the solution connected to the first selection steps. 
The used loss function in the architecture is the mean squared error (MSE), see Table~\ref{tab:technical}, comparing the solutions of the ordinary differential equations. 
Fig.~\ref{fig:MSE_selection} shows the MSE in the different selection steps for the two regarded data domains.

Fig.~\ref{fig:MSE_selection} shows comparable results for the two data settings. The MSE is decreasing with the selection steps even though this decrease is non-monotonous. 
While in the first selection step, many different mechanisms are involved, the last selection step includes only six parameters. 
In total, the MSE decreases to a lower value for the scenarios with more available data.

\section{Conclusion and further ideas}

The results show a successful model selection algorithm by using the longtime behavior of the solutions as qualitative data. 
The model selection algorithm reduces in each step the number of included mechanisms and selects in this way a model out of a model family with various mechanisms. 
The number of data points was rather small in the shown examples for imitating situations with expensive data in reality. 
As an effect of the limited data and rather small network structure, the approximation quality of the solutions of the learned model is not perfect and there are open questions to discuss. 
In this setting, the algorithm uses data that is rather qualitative than quantitative, and still provides the possibility of model selection. 

Further research ideas include the comparison of the loss function for classical sparse PINN approaches and the proposed model selection algorithm. Further, we plan to apply the approach to partial differential equations where the qualitative data is given by a characteristic solution behavior. In the partial differential equation setting, we want to include domain-decomposition techniques for separating dynamic and quasi-stationary time intervals or, in spatial heterogeneous settings, different reaction mechanisms in the spatial domains.

\begin{acknowledgement}
CR acknowledges funding for attending ENUMATH 2023 by the Seed Funding Programme of Technische Universit\"at Braunschweig, 2022 Interdisciplinary Collaboration: Strengthening Interdisciplinarity-- Expanding Research Collaboration.
\end{acknowledgement}
{\small \textbf{Competing Interests}
The authors have no conflicts of interest to declare that are relevant to the content of this chapter.}

  \bibliographystyle{acm} 
    \bibliography{references.bib}

\begin{thebibliography}{10}

\bibitem{Alber}
{\sc Alber, M., Buganza~Tepole, A., Cannon, W.~R., De, S., Dura-Bernal, S.,
  Garikipati, K., Karniadakis, G., Lytton, W.~W., Perdikaris, P., Petzold, L.,
  and Kuhl, E.}
\newblock Integrating machine learning and multiscale modeling—perspectives,
  challenges, and opportunities in the biological, biomedical, and behavioral
  sciences.
\newblock {\em npj Digital Medicine 2}, 1 (Nov. 2019), 115.

\bibitem{BruntonKutz-book}
{\sc Brunton, S.~L., and Kutz, J.~N.}
\newblock {\em Data-driven science and engineering: machine learning, dynamical
  systems, and control}.
\newblock Cambridge University Press, Cambridge New York, NY Port Melbourne New
  Delhi Singapore, 2019.

\bibitem{PySindy}
{\sc Brunton, S.~L., Proctor, J.~L., and Kutz, J.~N.}
\newblock Discovering governing equations from data by sparse identification of
  nonlinear dynamical systems.
\newblock {\em Proceedings of the National Academy of Sciences 113}, 15 (Apr.
  2016), 3932--3937.

\bibitem{SCIANN}
{\sc Haghighat, E., and Juanes, R.}
\newblock {SciANN}: {A} {Keras}/{TensorFlow} wrapper for scientific
  computations and physics-informed deep learning using artificial neural
  networks.
\newblock {\em Computer Methods in Applied Mechanics and Engineering 373\/}
  (Jan. 2021), 113552.

\bibitem{Meidani_2021}
{\sc Meidani, K., and Barati~Farimani, A.}
\newblock Data-driven identification of {2D} {Partial} {Differential}
  {Equations} using extracted physical features.
\newblock {\em Computer Methods in Applied Mechanics and Engineering 381\/}
  (Aug. 2021), 113831.

\bibitem{Rackauckas}
{\sc Rackauckas, C., Ma, Y., Martensen, J., Warner, C., Zubov, K., Supekar, R.,
  Skinner, D., Ramadhan, A., and Edelman, A.}
\newblock Universal {Differential} {Equations} for {Scientific} {Machine}
  {Learning}, Nov. 2021.
\newblock arXiv:2001.04385 [cs, math, q-bio, stat].

\bibitem{reisch-chemotactic}
{\sc Reisch, C., and Langemann, D.}
\newblock Chemotactic effects in reaction-diffusion equations for inflammation.
\newblock {\em Journal of Biological Physics 45}, 3 (Sept. 2019), 253--273.

\bibitem{reisch-chronification}
{\sc Reisch, C., and Langemann, D.}
\newblock Modeling the {Chronification} {Tendency} of {Liver} {Infections} as
  {Evolutionary} {Advantage}.
\newblock {\em Bulletin of Mathematical Biology 81}, 11 (Nov. 2019),
  4743--4760.

\bibitem{reisch-automative}
{\sc Reisch, C., and Langemann, D.}
\newblock Automative model selection and model certification for
  reaction-diffusion equations.
\newblock {\em IFAC-PapersOnLine 55}, 20 (2022), 73--78.

\bibitem{Schmiester_parametrization}
{\sc Schmiester, L., Weindl, D., and Hasenauer, J.}
\newblock Parameterization of mechanistic models from qualitative data using an
  efficient optimal scaling approach.
\newblock {\em Journal of Mathematical Biology 81}, 2 (Aug. 2020), 603--623.

\bibitem{Schmiester_qualitative}
{\sc Schmiester, L., Weindl, D., and Hasenauer, J.}
\newblock Efficient gradient-based parameter estimation for dynamic models
  using qualitative data.
\newblock {\em Bioinformatics 37}, 23 (Dec. 2021), 4493--4500.

\bibitem{Thomas2016}
{\sc Thomas, E., and Liang, T.~J.}
\newblock Experimental models of hepatitis {B} and {C} — new insights and
  progress.
\newblock {\em Nature Reviews Gastroenterology \& Hepatology 13}, 6 (June
  2016), 362--374.

\end{thebibliography}

\end{document}